\documentclass{amsart}

\usepackage{times, enumerate}
\usepackage[T1]{fontenc}

\usepackage[colorlinks=true]{hyperref}
\hypersetup{urlcolor=blue, citecolor=red}

\newcommand\Hom{\operatorname{Hom}}

\usepackage[normalem]{ulem}
\renewcommand\sout{\bgroup\markoverwith
{\textcolor{red}{\rule[0.7ex]{3pt}{1.4pt}}}\ULon}

\usepackage{times, amsmath, amsbsy, amssymb, verbatim,amscd}
\usepackage{color,graphicx,fancybox}
\DeclareGraphicsExtensions{.jpg,.pdf,.png,.ps,.eps}

\RequirePackage{color}
\definecolor{DarkBlue}{rgb}{0.0,0.0,0.55}
\definecolor{LightYellow}{rgb}{1.0,1.0,0.8}
\definecolor{mygreen}{rgb}{0.0,0.501,0.0}
\definecolor{mybrown}{rgb}{0.647,0.164,0.164}
\definecolor{lightsteelblue3}{rgb}{0.634,0.710,0.804}
\definecolor{gray}{rgb}{0.2,0.2,0.2}

\newcommand\CC{\mathbb C}

\newcommand\RR{\mathbb R}

\newcommand\ZZ{\mathbb Z}

\newcommand{\maC}{\mathcal C}

\newcommand{\maG}{\mathcal G}

\newcommand{\maK}{\mathcal K}

\newcommand{\maV}{\mathcal V}

\newcommand{\pa}{\partial}
\newcommand\ede{ \, := \, }
\newcommand\seq{ \, = \, }

\newcommand\oM{\overline{M}}

\newcommand{\CI}{\mathcal C^{\infty}}
\newcommand\CIc{\mathcal{C}^\infty_c}

\newcommand\m[1]{$#1$}

\newcommand\Diff{\operatorname{Diff}}

\newtheorem{theorem}{Theorem}[section]

\newtheorem{problem}[theorem]{Problem}

\theoremstyle{definition}

\newtheorem{definition}[theorem]{Definition}

\newtheorem{remark}[theorem]{Remark}

\title[Non-compact manifolds]{Analysis on noncompact manifolds and Index Theory:
Fredholm conditions and Pseudodifferential operators}

\author[I. Beschastnyi]{Ivan Beschastnyi} \address{
  CIDMA, University of Aveiro, Campus Santiago, 3810-193, Aveiro, Portugal}
\email{i.beschastnyi@gmail.com}

\author[C. Carvalho]{Catarina Carvalho} \address{{CEAFEL, Dep. Matem\'{a}tica,
  Instituto Superior T\'{e}cnico, University of Lisbon, Av. Rovisco
  Pais, 1049-001 Lisbon, Portugal 
}}
\email{catarinaccarvalho@tecnico.ulisboa.pt}

\author[V. Nistor]{Victor Nistor} \address{
Universit\'{e} de Lorraine, CNRS, IECL, F-57000 Metz, France}
\email{victor.nistor@univ-lorraine.fr}

\author[Y. Qiao]{Yu Qiao} \address{School of Mathematics and
  Statistics, Shaanxi Normal University, Xi'an, 710119,
China} \email{yqiao@snnu.edu.cn}

\thanks{I.B and C.C were supported by FCT--Funda\c c\~ao para a Ci\^encia e Tecnologia, 
under projects UIDB/04106/2020 and UIDB/04721/2020, respectively. V.N. was supported by ANR grant 
OpART and Y.Q was supported by NSFC (11971282).}

\begin{document}

\begin{abstract} 
  We provide Fredholm conditions for compatible differential 
  operators on certain Lie manifolds (that is, on certain
  possibly non-compact manifolds with nice ends). 
  We discuss in more detail the case of manifolds with cylindrical, 
  hyperbolic, and Euclidean ends, which are all covered by 
  particular instances of our results. We also discuss applications
  to Schr\"odinger operators with singularities of the form 
  $r^{-2\gamma}$, $\gamma \in \RR_+$.
\end{abstract}

\maketitle


\section{Introduction}
We review in this note an approach to \emph{Analysis on non-compact manifolds}
that is based on \emph{Lie algebras of vector fields} and \emph{compactifications of manifolds.} 
This note is a revised and enhanced version of the talk delivered by V.N. at the ``Methusalem 
Colloquium and Geometric Analysis Seminar,'' at Ghent University, in June 2022. 
The mathematical results included here are based, to a large extent, on the paper
\cite{CNQ}, but they are also enhanced by recent results from \cite{BCNQ}. Among the 
questions that we will touch upon are: 
\begin{itemize}
  \item the invertibility of differential operators 
  (Hadamard well-posedness), 
  \item the regularity of solutions of elliptic PDEs, 
  \item the construction of the algebras of pseudodifferential
  operators containing the inverses of the elliptic operators in question, and, most notably, 
  \item the \emph{Fredholm alternative,} the topic around which we organize our paper. 
\end{itemize}
\emph{Informally stated,} our main result on the Fredholm alternative is:

\begin{theorem}\label{thm.informal}
On a ``nice'' manifold, an adapted (pseudo)differential
operator is \emph{Fredholm} if, and only if, 
\begin{enumerate}[(i)]
    \item it is \emph{elliptic} and
    \item \emph{all its ``limit operators'' are invertible.}
\end{enumerate}
\end{theorem}

A more formal statement will be included below (Theorem \ref{thm.main.Fredholm}).
In the above result, by ``elliptic'' we mean a pseudodifferential operator whose 
principal symbol is invertible outside of the zero section. On a compact manifold,
there are no ``limit operators,'' and the Fredholm property reduces to ellipticity.
In the non-compact case, the 
\emph{limit operators do appear and play an important role.} The appearance of 
limit operators is thus one of the main differences between the analysis on 
compact and on non-compact manifolds. A large part of this note is devoted to 
presenting a more formal version of the above theorem (Theorem \ref{thm.informal}), 
as well as to explaining what the statement of that theorem becomes for some standard 
classes and then for some new classes
of manifolds. We consider in this regard manifolds that have ends that are asymptotically 
of one of the following forms: 

  \ (i) cylindrical, 
  \ (ii) hyperbolic, or 
  \ (iii) Euclidean (more general: conical).\\
In particular, we pay special attention to the description 
and properties of the limit operators for these manifolds. 
The new type of examples is related to Schr\"odinger operators 
with potentials with non-integer order singularities at the origin. These examples
rely on a generalization of the setting of \cite{CNQ} that leads to more general 
pseudodifferential calculi \cite{BCNQ}.

There are many reasons to study Fredholm operators, one being that they are useful
also in obtaining \emph{well-posedness results}, for instance, by the following 
paradigm: to prove that an operator $P$ is invertible, it suffices to prove that it
is Fredholm, that it has index zero, and that it has vanishing kernel. 
An example where this procedure was 
used is in the paper by M. Mitrea and M. Taylor \cite{MitreaTaylorJFA} studying
the layer potentials on smooth subdomains of compact manifolds. Our results can thus
be used to extend the work of Mitrea and Taylor to some non-compact manifolds \cite{GKN1}. 

Let us stress that, although we are advertising mostly smooth (non-compact) 
manifolds, the results presented here can be useful also for (pseudo)differential operators on singular 
spaces, since one can treat singular spaces by doing analysis on their smooth part, which is 
regarded as a non-compact manifold, possibly with a different metric, usually conformally 
equivalent to the original one. This is the case with \emph{the passage from
a manifold with conical points to a manifold with cylindrical ends,} which was classically
done via the  so called ``Kondratiev transform'' $r = e^{t}$ and is explained in 
the text. 

Compared to the original talk, we have removed some repetitions and 
we have included some additional details, examples, and results; we have tried, however, to maintain
as much as possible the style (and idiosyncrasies) of the original talk. In particular,
the plan of the presentation is to successively discuss the Fredholm alternative
and related topics for:
\begin{enumerate}[(i)]

  \item \emph{Smooth, compact manifolds} (the classical case).
  \item \emph{Manifolds with cylindrical ends} (an almost classical case).
  The main examples here are the Laplacian in polar and generalized spherical coordinates.
  \item \emph{Other classes of manifolds}, for which we stress the similarities and differences to manifolds with 
  cylindrical ends. The main examples here are the Laplacian in cylindrical coordinates
  and in flat, Euclidean
 coordinates.
\end{enumerate}
As in the original talk (in any talk, for that matter), it was not realistic 
to include a complete list of references. This is unfortunate, because very many people 
have worked on analysis on non-compact manifolds. Since we are not including enough
references, let us stress that, \emph{unless explicitly stated otherwise,} none of the 
results below belong to us. We thank Cipriana Anghel, Sergiu Moroianu, Elmar Schrohe, 
and J\"org Seiler for useful discussions.

\section{Motivation and some classical results}

Let us see what Theorem \ref{thm.informal}
becomes in some classically well-understood cases, namely those of compact manifolds and of manifolds
with cylindrical ends.

\subsection{A classical case: smooth, compact manifolds}
For pedagogical reasons, it will be useful to recall first the 
Fredholmness result on smooth, compact manifolds without boundary. For simplicity, we will 
assume from now on that all our manifolds are \emph{smooth and complete Riemannian}
(except the manifolds with conical points).
Thus we can define the Sobolev spaces on these manifolds using the powers of the Laplacian. 
Let $P$ be an order $m$, classical, pseudodifferential operator 
on a compact manifold \m{M}. Then $P : H^s(M) 
\to H^{s-m}(M)$ is bounded. In general, we will consider operators acting 
between sections of vector bundles $E$ and $F$. 
The following result is classical:

\begin{theorem}\label{thm.classical}
    Assume that \m{M} is \emph{smooth, compact, without boundary,} and that $P$ is an order $m$, 
    classical, pseudodifferential operator. Then
    \m{P : H^s(M; E) \to H^{s-m}(M; F)} is \emph{Fredholm} if, and only if, it is \emph{elliptic.}
\end{theorem}

We stress again that, in the compact case, \emph{there are no limit operators.}
As we will see in the next subsection, this is \emph{not} the case in the non-compact 
case.

\subsection{Motivating example: cylindrical ends}
The simplest, motivating example of a non-compact manifold is that of a 
manifold with \emph{cylindrical ends.} We begin with three examples before giving 
the formal definition.

\subsubsection{Polar coordinates} The Laplacian on \m{\RR^2} in 
\emph{polar coordinates} \m{(r, \theta)} is 
\begin{equation*}
  \Delta_{\RR^2} u \seq  r ^{-2} \big(( r  \pa_{ r })^2 u +
    \pa_\theta^2u \big). 
\end{equation*}
Ignoring \m{ r^{-2},} we obtain the differential operator $(r  \pa_{r})^2  +
    \pa_\theta^2$
acting on \m{[0, \infty) \times S^{1} \ni (r, \theta),} which is a  
\emph{degenerate elliptic} differential operator acting on a manifold with boundary. We follow 
Melrose's observation \cite{MelroseAPS} that this operator is 
generated by the vector fields $$r \pa_r \ \mbox{ and } \ \pa_\theta\,,$$
which are tangent to the boundary $\{0\} \times S^1$, 
an observation that will be useful in generalizations. 

This example can be treated with the help of the \emph{Kondratiev transform}
$t = \log r,$ which transforms $r \pa_ r$ into $\pa_t.$ The operator 
$(r  \pa_{r})^2  + \pa_\theta^2$ then becomes $\pa_{t}^2  + \pa_\theta^2$. 
Let $S^{n-1}$ be the \emph{unit sphere} in \m{\RR^n.}
The Kondratiev transform maps the domain $(0, \infty) \times S^{1}$
to $\RR \times S^{1}$. The advantage of the Kondratiev transform is that
the resulting operator \m{\pa_{t}^2  + \pa_\theta^2} is nothing but \m{\Delta_{\RR \times S^1},} 
the \emph{ Laplacian} for the metric \m{(dt)^2 + (d\theta)^2} on $\RR \times S^1$. 
The Fourier transform in \m{t} maps \m{\pa_t} to \m{i \tau} and hence it maps 
\m{\pa_{t}^2 + \pa_\theta^2} to \m{ - \tau^2 + \pa_\theta^2.} In turn, this operator
can be understood via the spectral theory of \m{\pa_\theta^2} on \m{S^1}, 
which is well known. Anticipating, $\RR \times S^1$ is the simplest non-trivial example of 
a ``manifold with cylindrical ends.'' For general manifolds with cylindrical ends, the 
Fourier transform will be applied not to the operator itself, but rather to its limit
operators, which are also $\RR$-invariant operators.

\subsubsection{The Black-Scholes operator}
An example that can be treated similarly is that of
\begin{equation*}
    \pa_t u \, + \, \, \frac{\sigma^2}{2} (x\pa_x)^2 u + \Big( r - \frac{\sigma^2}{2} \Big) 
    x \pa_x u - ru \,,
\end{equation*}
the \emph{Black-Scholes operator}. 
This example is, in a certain sense, even simpler than that of the Laplace operator
in polar coordinates, since it is an example in one dimension.

\subsubsection{Generalized spherical coordinates} 
A related example is the \emph{Schr\"{o}dinger operator} 
$\Delta_{\RR^n} + Z/\rho$ on \m{\RR^n,}
$Z \in \CC$, in {\em (generalized) spherical coordinates}
\m{(\rho, x') \in (0, \infty) \times S^{n-1}}:
\begin{equation*}
  -\Big(\Delta_{\RR^n} + Z\rho^{-1} \Big) u \seq -\rho^{-2} \big(
    (\rho \pa_{\rho})^2 u + (n-2)\rho \pa_{\rho} u +
    \Delta_{S^{n-1}}u + Z \rho u \big).
\end{equation*}
In this example, the operator is considered on $[0,\infty) \times S^{n-1}$
and is generated by $\rho \pa_\rho$ and vector fields independent of $\rho$
(which yield the Laplacian on the sphere $S^{n-1}$). 
The relevant geometry for this operator is again that 
of manifolds with cylindrical ends, since it can be generated by vector 
fields that are tangent to the boundary $\{0\} \times S^{n-1}$
of $[0,\infty) \times S^{n-1}$. Note that, after ignoring the factor $\rho^{-2}$,
the resulting operator no longer has a singular potential. This example can also be treated 
via the Kondratiev and Fourier transforms.

\subsubsection{Manifolds with cylindrical ends} Let us now describe the general 
setting encompassing the previous three examples. 
Let \m{\oM} be a \emph{smooth, compact manifold with boundary} \m{\pa \oM} and 
let $r \ge 0$ be a smooth function on 
$\oM$ such that $\pa \oM = r^{-1}(0)$ and $dr \neq 0$ on $\pa \oM$. Up to Lipschitz equivalence,
a \emph{manifold with cylindrical ends} is one that is isometric to 
$M := \oM \smallsetminus \pa \oM$ endowed with a metric that, 
near the boundary, is of the form 
\begin{equation}\label{eq.def.gcyl}
    g_{cyl} \ede \frac{dr^2}{r^2} + h  \,,
\end{equation}
where $h$ is a semi-definite tensor that restricts to a true metric on $\pa \oM$.
In local coordinates $x = (r, x_2, \ldots, x_n)$ near the boundary, 
the Sobolev spaces $H^m(M; g_{cyl})$ associated to the metric $g_{cyl}$
identify with the Babu\v{s}ka-Kondratiev (or weighted Sobolev) spaces 
\begin{equation} 
  \maK_{a}^{m}(M) \ede \{r^{\alpha_1 - a} \pa^{\alpha} u
  \in L^2(M),\, |\alpha| \le m\} 
\end{equation}
 with $a = \dim(\oM)/2$. These function spaces (for all $a$) arise naturally if we study the manifold with 
conical points $\oM/\sim$, where $\sim$ collapses each component of $\pa \oM$ to a point. 

Let us now turn to the definition of \emph{differential operators.} 
Assume, first that $\oM = [0, \infty) \times \RR^{n} \ni (r, x')$ and let 
\m{\alpha = (\alpha_1, \alpha_2, \ldots, \alpha_{n}) = (\alpha_1, \alpha') \in \ZZ_+^{n+1}} be a 
generic multi-index. Then the ``right differential operators'' for the 
cylindrical ends geometry are the \emph{totally characteristic differential operators} $P$
of the form
\begin{equation}\label{eq.tot.char} 
    P \seq \sum_{|\alpha| \le m} a_{\alpha}(r, x') (r \pa_r)^{\alpha_1} \pa_{x'}^{\alpha'}\,.
\end{equation}
Let us endow $M$ with the metric $\frac{(dr)^2}{r^2} + (dx')^2$. 
Then, the Laplacian for this metric is a totally characteristic differential operator. 
Furthermore, after the Kondratiev transform $r = e^t$, 
this metric becomes just the usual Euclidean metric 
on $\RR^{n+1} \ni (t, x')$. On the other hand, the differential operator $P$ becomes
\begin{equation}
    P \seq \sum_{|\alpha| \le m}  a_{\alpha}(e^t, x') \pa_t^{\alpha_1} \pa_{x'}^{\alpha'}\,,
\end{equation}
which has the property that its coefficients have limits $a_{\alpha}(0, x')$
as $t \to - \infty$. This leads us to an 
important new ingredient, namely the ``$b$-normal'' (simply, ``normal'') operator
\begin{equation}
    \tilde P \ede \sum_{|\alpha| \le m}  a_{\alpha}(0, x') \pa_t^{\alpha_1} \pa_{x'}^{\alpha'}\,,
\end{equation}
obtained by ``freezing the coefficients'' at $t = - \infty$.
The normal operator is an important instance of a ``limit operator'' mentioned 
earlier. (Notice that $r \pa_r$ was replaced with $\pa_t$, so we have used 
the Kondratiev transform implicitly.) The Fourier transform in $t$ then yields the 
\emph{indicial family} 
$\widehat P(\tau) := \sum_\alpha a_{\alpha}(0, x') (i \tau)^{\alpha_1} \pa_{x'}^{\alpha'}\,.$
Similar constructions ($\tilde P$ and $\widehat P$) can be defined on any manifold 
with cylindrical ends by localization, and we obtain a normal operator for each 
connected component $L$ of the boundary $\pa \oM$ of $\oM$.  We have the following (almost classical) result
\cite{Kondratiev67, MelroseMendoza, SchulzeBook91}:

\begin{theorem}                               
Let $D$ be an order $m$, totally characteristic differential operator on 
$\oM$. The map 
$D : H^{s}(M; E, g_{cyl}) \to H^{s-m}(M; F, g_{cyl})$ is Fredholm if, and 
only if, 
\begin{enumerate}[(i)]
\item it is elliptic and
\item all its normal operators are invertible $H^{s}(L; E) \to H^{s-m}(L; F)$.
\end{enumerate} 
\end{theorem}

An approach to the analysis on manifolds with cylindrical ends is via the $b$-calculus
of Melrose \cite{MelroseAPS} and Schulze \cite{SchulzeBook91}. See also \cite{LauterSeiler}.

\begin{figure}
    \includegraphics[width=1.55in]{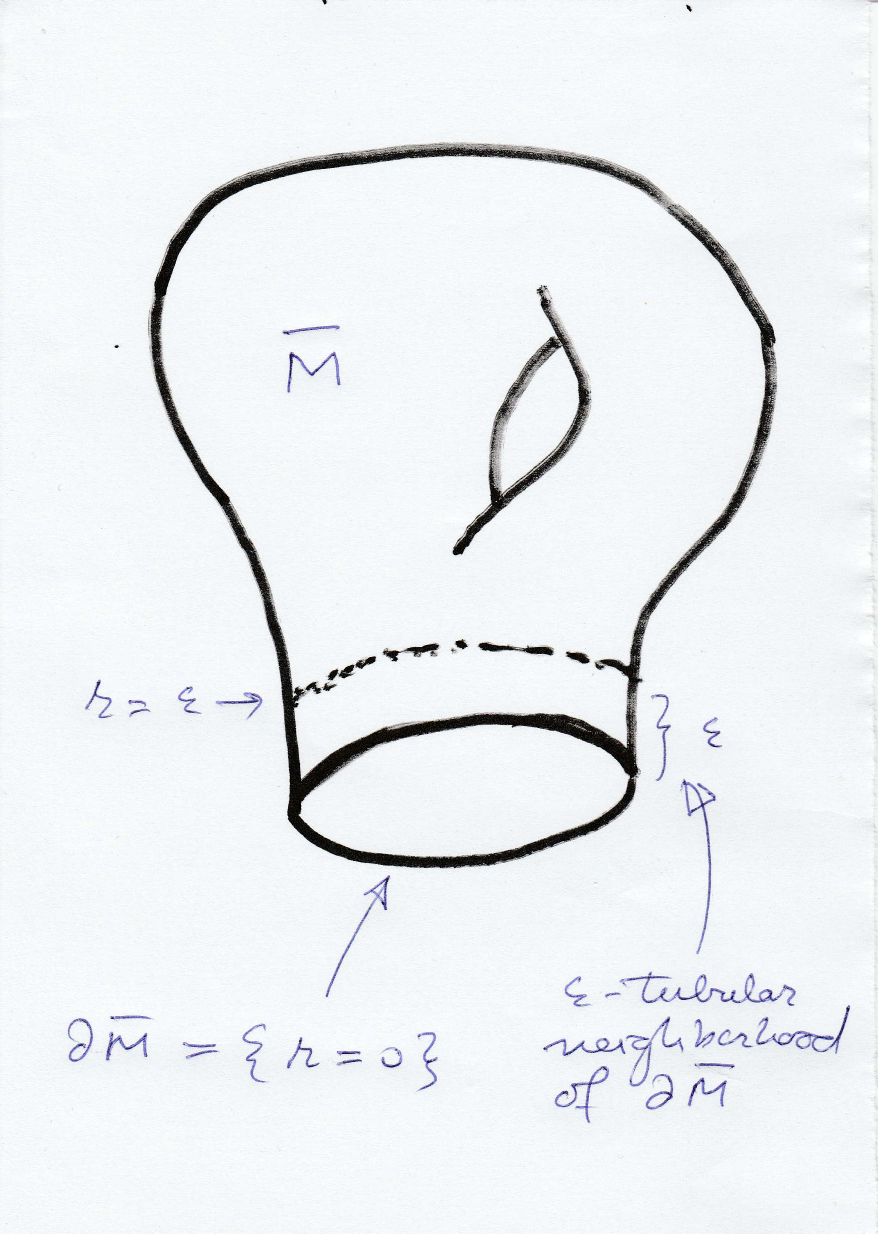}  \qquad \
    \includegraphics[width=1.55in]{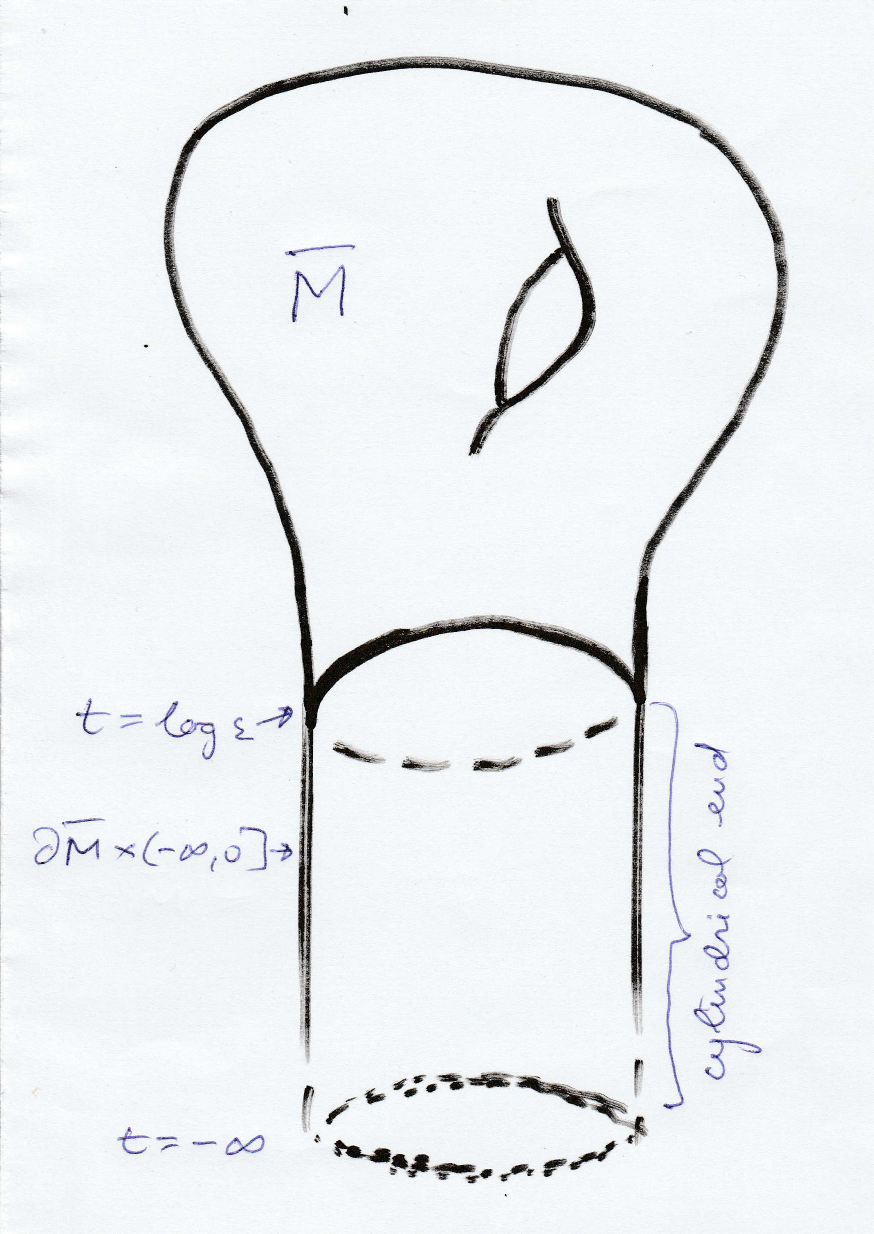}
    \caption{The construction of a manifold with cylindrical ends.}
\end{figure}

\section{Fredholm operators on manifolds with nice ends (Lie manifolds)}

The Laplace operator in cylindrical coordinates in $\RR^3$ or the Schr\"odinger 
operators with singular potentials of the form $Z\rho^{-2\gamma}$, $\gamma \in 
\RR_+ \smallsetminus \{0, 1/2, 1\}$, do not fit into the framework of manifolds with cylindrical ends.
This raises the question of defining compatible pseudodifferential calculi and 
obtaining the associated Fredholm alternative (or Fredholm conditions) on more 
general manifolds.
It does not seem possible to obtain convenient characterizations of 
Fredholm operators on general non-compact manifolds. We will thus restrict ourselves
to a class of non-compact manifolds, which we will call ``nice manifolds'' for
the purpose of this paper. This class of manifolds consists of manifolds with 
nice ends and are, for the most part,
``Lie manifolds,'' a class of manifolds that 
we introduce next. While we do not define precisely which manifolds are nice 
(this is rather technical, see \cite{BCNQ, CNQ}), we do discuss several examples,
among which the ones mentioned in the introduction, namely manifolds with (asymptotically): cylindrical ends,
hyperbolic ends, and Euclidean
 ends.

\subsection{Lie manifolds}
We now introduce and discuss Lie manifolds, their geometry, and, most
importantly, their \emph{associated differential operators.}

\subsubsection{Definition of Lie manifolds}
Assume that we are given a compact manifold with corners $\oM$, whose 
interior is $M := \oM \smallsetminus \pa \oM$, and a subspace
\begin{equation}\label{eq.def.maV1}
    \maV \subset \maV_b(\oM) \ede \{ X \in \CI(\oM; T\oM) \mid X 
    \mbox{ \emph{tangent} to all faces of } \oM \} \,.
\end{equation}

\begin{definition}[\cite{aln1}] 
The pair \m{(\oM, \maV)} is a Lie manifold if: 
\begin{enumerate}[(i)]

\item $\maV \subset \maV_b(\oM)$ is closed under the Lie bracket $[\;,\;]$;

\item $\maV$ is a finitely-generated, projective \m{\CI(\oM)}--module;

\item $\CIc(M; TM) \subset \maV$ (recall that $M := \oM \smallsetminus \pa \oM$).
\end{enumerate}
\end{definition}

This definition is based on earlier, similar constructions due to Connes,
Cordes, Kondratiev, Mazzeo, Mazya, Melrose, Plamenevskij (a relevant old 
article with Mazya is described in \cite{NP}), Schrohe, Schulze, 
Skandalis, and many others. 
Informally, a {\em Lie manifold} \m{(\oM, \maV)} is a manifold \m{M} with:
\begin{enumerate}[(i)]
    \item a \emph{compactification} \m{\oM} such that 
    \m{M = \oM \smallsetminus \pa \oM,} whose role is to control 
    the behavior at infinity of the coefficients of our differential operators, and 

    \item a \emph{Lie algebra} of vector fields \m{\maV \subset \maV_b} on this
    compactification, whose role is to define the ``nice'' differential operators.
\end{enumerate}

\begin{remark}
Let \m{(\oM, \maV)} be a Lie manifold. Saying that \m{\maV} is a {\em projective} 
\m{\CI(\oM)}-module means that it is stable under multiplication by functions in 
\m{\CI(\oM)} and it has a basis around each point of \m{\oM.} In particular, $\maV$ is a 
complex vector space. Saying that \m{\CIc(M; TM) \subset \maV} means that there are 
no ``obstructions in the interior'' for the vector fields in $\maV$. Compact (smooth) manifolds
and manifolds with cylindrical ends are Lie manifolds. These and other examples will be
discussed below in Section \ref{ssec.examples}.
\end{remark}

\subsubsection{Lie manifolds and geometry}
As in \cite{aln1}, since $\maV$ is a projective $\CI(\oM)$-module, the 
Serre-Swan theorem gives the existence of a vector bundle \m{A\to \oM}
such that
\begin{equation}\label{eq.def.A}
  \maV \simeq \Gamma(\oM; A)\,.
\end{equation}
The Lie algebra 
structure on the sections of $A$ means that it is a \emph{Lie algebroid}.
The inclusion $\CIc(M; TM) \subset \maV$ means that \m{A = TM} in the interior 
$M := \oM\smallsetminus \pa \oM$ of $\oM$.
Thus, once we have chosen a metric on \m{A}, that metric will induce a metric on \m{TM}
(or, which is the same thing, on \m{M}). A metric obtained in this way will be 
called \emph{compatible}, and is unique up to 
Lipschitz equivalence. For a smooth bundle $E \to \oM$, 
we define the Sobolev spaces $H^s(M; E)$ using the powers of the Laplacian for any compatible metric.

\subsubsection{Differential operators}Our main interest lies in the algebra
\begin{equation*}
    \Diff_{\maV}(\oM),
\end{equation*} 
which consists of the differential operators generated by \m{\CI(\oM)} and by 
derivatives in \m{\maV}. Given two vector bundles $E, F \to \oM$,
we have a similar definition of $\Diff_{\maV}(\oM; E, F)$, 
the set of differential operators 
generated by $\CI(\oM; \Hom(E; F))$ and $\maV$ acting on sections 
of $E$ with values sections of $F$. 
All the geometric operators on $M$ with a compatible metric will belong to some space of the 
form $\Diff_{\maV}(\oM; E, F)$. For instance, the Levi-Civita connection on $M$ extends 
to a differential operator in $\Diff_{\maV}(\oM; A, A^* \otimes A)$ (an $A$-connection).
Hence, since $\oM$ is compact, its interior $M$ will have bounded curvature (together with its 
covariant derivatives). 

\subsection{Formulation of the main result: Fredholm conditions on ``nice'' 
manifolds} We are now closer to giving a precise statement for the 
Fredholm conditions in theorem \ref{thm.informal}. 
Let \m{(\oM, \maV)} be a Lie manifold and $D\in \Diff_{\maV}(\oM; E, F)$. 
If $D$ has order $m$, then it defines a continuous map $D : H^s(M; E) \to H^{s-m}(M; F)$. 
To \m{(\oM, \maV)} and \m{D} we associate: 
\begin{enumerate}[(i)]
	\item spaces $Z_\alpha$, $\alpha \in I$, 
	the \emph{orbits} of $\maV$ on $\oM$ (they do not dependent on $D$);
	\item isotropy groups $G_\alpha$, $\alpha \in I$,
	(also independent of \m{D}); and
 	\item \m{G_\alpha}-invariant differential operator \m{D_\alpha}
          on \m{Z_\alpha \times G_\alpha}
\end{enumerate}
as follows. Let ${\maV_x} := \{ X \in \maV\vert\, X(x) = 0\}$ and 
${I_x} := \{ f \in \CI(\oM)\vert\, f(x) = 0\},$ for $x\in \pa \oM$.  Then $\maV_x /I_x \maV$
is a Lie algebra and, for $x \in Z_\alpha$, $G_\alpha$ is the simply-connected Lie group that integrates it:
\begin{equation*} 
    Lie (G_\alpha) \seq {\maV_x /I_x \maV}\,.
\end{equation*}
For each $\alpha \in I$, the operator $D_{\alpha}$ is obtained by restricting $D$ to the orbit 
$Z_\alpha$ and letting the vector fields act on $G_\alpha$ via the 
map $\maV_x \to Lie(G_\alpha)$. The operators \m{D_\alpha} are the 
\emph{limit operators}.
We can now formulate our main result \cite{BCNQ, CNQ}

\begin{theorem} \label{thm.main.Fredholm}
  Let $D \in  \Diff_{\maV}(\oM; E, F)$ be of order $m$ and assume that $(\oM, \maV)$
  is ``nice.'' We have that $D : H^s(M; E) \to H^{s-m}(M; F)$ is Fredholm
  if, and only if, 
  \m{D} is elliptic and all \m{D_{\alpha}}, \m{\alpha \in I}, are invertible.
\end{theorem}

Our Lie manifold $(M, \maV)$ is \emph{nice} if there exists a \emph{Hausdorff} 
Lie groupoid $\maG$ such that
\begin{itemize}
    \item its Lie algebroid $A(\maG) \simeq A$ of Equation \eqref{eq.def.A}
    (so $\Gamma(\oM; A) = \maV$) and 
    \item $\maG$ satisfies the Effros-Hahn conjecture (a statement about 
    $C^*(\maG)$ that implies that the regular representations
    of $C^*(\maG)$ determine its invertible elements).
\end{itemize}
The proof is obtained by considering the norm closure $\overline{\Psi}^0(\maG)$
of the algebra of order zero pseudodifferential operators on $\maG$ \cite{BCNQ, CNQ}. 
Let $\maK$ be the algebra of compact operators on $L^2(M)$. Recall that an operator 
is Fredholm if, and only if, it is invertible modulo $\maK$. We thus want to 
characterize the invertible elements of $\overline{\Psi}^0(\maG)/\maK$. Since 
$\overline{\Psi}^0(\maG)/C^*(\maG) \simeq \maC(S^*A)$ via the principal symbol
and $\maG$ satisfies the the Effros-Hahn conjecture, in addition to 
the principal symbol, it is enough to look at the regular representations of 
$C^*(\maG)/\maK$. Each of these regular representations yields a limit operator.
This completes the proof. 
A useful property here is that 
nice manifolds are closed under general blow-ups \cite{BCNQ, CNQ}.

\begin{figure}
    \includegraphics[scale=.065]{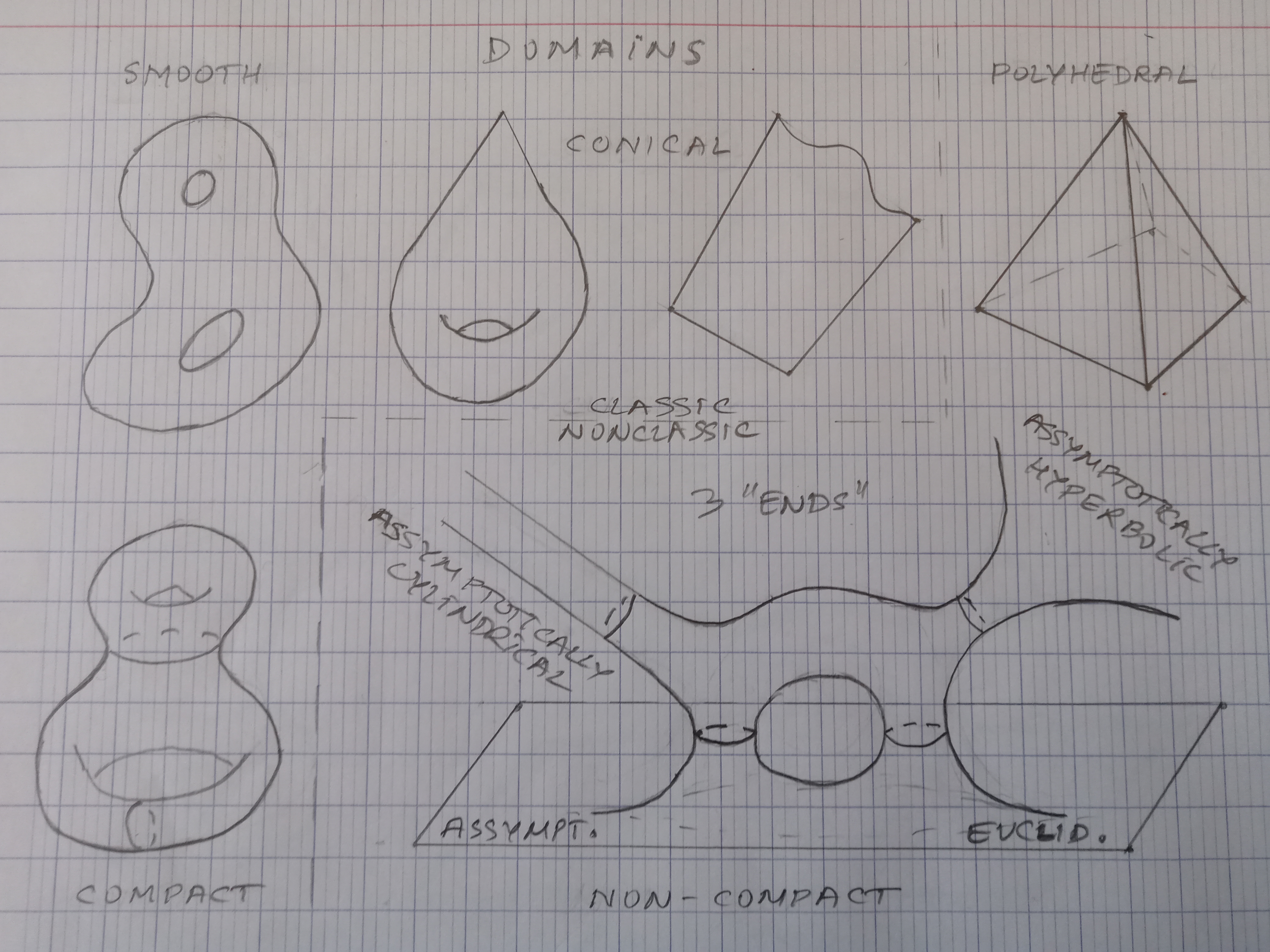}
    \caption{Various types of singular spaces or domains that can be treated
    using blow-ups and Lie manifolds and an example of a non-compact manifold that 
    has three types of ends: cylindrical, hyperbolic, and Euclidean.}
\end{figure}

\subsection{Examples and applications}\label{ssec.examples}
The example ``zero'' is that of a smooth, compact manifold without boundary $M$.
Then $M$ is a Lie manifold with $\oM = M$ and 
$\maV = \Gamma(M; TM)$ (all smooth vector fields on $M$). In particular, we have 
$A = TM$. The boundary $\pa \oM := \oM \smallsetminus M$ is empty, which is consistent 
with the fact that there are no limit operators. Then Theorem \ref{thm.main.Fredholm} 
becomes simply the (well known) Theorem \ref{thm.classical}. 
Thus the example of compact 
manifolds, while it fits our setting, is too simple to illustrate the general theory, 
and hence we will concentrate next on other, more relevant examples.
In all our following examples, \m{\oM} will be a smooth, compact manifold with boundary $\pa \oM$
and we shall concentrate on a tubular neighborhood \m{U = [0, \epsilon) \times \pa \oM \ni ({r}, y)}
of the boundary, because the interior $M := \oM \smallsetminus \pa \oM$ can be 
treated with classical tools. Here are 
our main three examples:

\begin{enumerate}[(i)]
  \item 
  Let $\maV := \maV_b$, the set of smooth vector fields on $\oM$ that are
  \emph{tangent} to the boundary $\pa \oM$ of $\oM$, as before. 
  We have $\maV_b = \CI(\oM){r \pa_r} + \sum \CI(\oM){\pa_y}$ on $U$, so the 
  projectivity condition is satisfied and we obtain a Lie manifold $(\oM, \maV_b)$.
  Any metric of the form $g = r^{-2}(dr)^2 + h$ as in Equation \eqref{eq.def.gcyl}
  will be compatible (with the Lie manifold structure) on $M$, and hence 
  $M$ is a \emph{manifold with cylindrical ends.}
  The bundle \m{A = {}^{b}TM} was considered by Melrose. The index set 
  \m{I} consists of the connected components of $\pa \oM$ and \m{G_\alpha = \RR}. The differential 
  operators $\Diff_{\maV}(\oM)$ are the totally characteristic differential operators on $\oM$ 
  (already discussed) and the limit operators $D_\alpha$ are the normal operators of $D$ defined 
  earlier.
 
  \item Let next $\maV := \maV_0 := r \CI(\oM; T\oM)$, the set of smooth vector fields on $\oM$ that 
  \emph{vanish} at the boundary $\pa \oM$ of $\oM$.  Then $\maV_0 = \CI(\oM){r \pa_r} + \sum \CI(\oM){r\pa_y} $ 
  on $U$, so the projectivity condition is satisfied and, again, 
  $(\oM, \maV_0)$ is a Lie manifold. The choice of a metric on $A$
  (a compatible metric on $M$) makes the interior $M := \oM \smallsetminus \pa \oM$ of $\oM$
  a so called {\em asymptotically hyperbolic manifold}, since the metric is of the form 
  $g = \frac{h}{{r^2}}$ with $h$ a true metric on $\oM$. The orbits $Z_\alpha \subset \pa \oM$ 
  are reduced to points (so $I = \pa \oM$)  and $G_\alpha = T_\alpha \pa \oM {\rtimes} \RR$, 
  which is a \emph{non-commutative} 
  group (obtained as the semi-direct product from the action of $\RR$ on $T_\alpha \pa \oM$ by dilations).
  
  \item To round up our podium, let $\maV \ede \maV_{sc} \ede r \maV_b \subset \maV_0$.
  From example (i), we see that $\maV_{sc} \ede \CI(\oM){r^{2} \pa_r} 
  + \sum \CI(\oM){r\pa_y}$ on $U$, so the projectivity condition is yet again 
  satisfied to yield a Lie manifold $(\oM, \maV_{sc})$. Euclidean, asymptotically
  Euclidean, and asymptotically conical spaces are modeled by this type of 
  Lie manifolds. The orbits \m{Z_\alpha = \{\alpha\}} are again reduced to points
  (so again $I =  \pa \oM$).
  However, this time, $G_\alpha = T_\alpha \pa \oM {\times} \RR$ is \emph{commutative.}
\end{enumerate}

The Laplacian $\Delta_{\RR^3} u = r^{-2} \big[(r \pa_{r})^2 u 
+ \pa_\theta^2u  + (r \pa_z)^2  u \big]$ in {\em cylindrical coordinates} $(r, 
\theta, z) \in [0, \infty) \times S^1 \times \RR$ is closely related 
to the example (ii) above. Ignoring the factor $r^{-2}$,
the relevant operator is generated by the vector fields $r \pa_r,$ $\pa_\theta,$ and $r \pa_z,$
which are again tangent to the boundary. 
This example was one of the original motivations to considering a framework more general than 
that of manifolds with cylindrical ends.

We conclude this section with an example that goes even beyond Lie manifolds.
Namely, let us we consider the Schr\"odinger operator with potential $V = \rho^{-2\gamma}V_0$,
with $V_0$ \emph{smooth in generalized spherical coordinates} $(\rho, x')$, $\rho \in [0, \infty)$, 
$x' \in S^{n-1}$, and $\gamma \ge 0$:
\begin{equation}\label{eq.schr1}
  \Delta + V \seq \rho^{-2}\big[ (\rho \pa_\rho)^2 
  + (n-2) \rho \pa_\rho + \Delta_{S^{n-1}} +  \rho^{2-2\gamma}V_0\big]
\end{equation} 
Let us assume first 
that $2\gamma \in \{0, 1, 2\}$. Ignoring the factor $\rho^{-2}$, the resulting 
operator 
\begin{equation}\label{eq.schr2}
    (\rho \pa_\rho)^2 + (n-2) \rho \pa_\rho + \Delta_{S^{n-1}} +  \rho^{2-2\gamma}V_0
\end{equation}
is in the $b$-calculus \cite{MelroseAPS, SchulzeBook91}. This property can be used, among
other applications, to study 
the domain of $\Delta + V$ and the regularity of 
its eigenfunctions. If $2\gamma \in \RR \smallsetminus \{0, 1, 2\}$, the operator
$\Delta + V$ can still be studied with a modified pseudodifferential calculus \cite{BCNQ}. 
 For instance, if $\gamma > 1$ and if we write
  $\Delta + V = \rho^{-2\gamma}\big[ (\rho^{\gamma} \pa_\rho)^2 
  + (n-1-\gamma) \rho^{\gamma-1} \rho^\gamma \pa_\rho + \rho^{2\gamma-2}\Delta_{S^{n-1}} + V_0\big],$
then the resulting operator 
\begin{equation}\label{eq.schr4}
  (\rho^{\gamma} \pa_\rho)^2 
  + (n-1-\gamma) \rho^{\gamma-1} \rho^\gamma \pa_\rho + \rho^{2\gamma-2}\Delta_{S^{n-1}} + V_0
\end{equation}
is in the $c_{\gamma, \gamma-1}$-calculus of \cite{BCNQ} (recall that $V_0$ is smooth in polar coordinates).

Other examples come from semi-Riemannian and sub-Riemannian geometry \cite{B_RIMA}.

\section{Pseudodifferential operators and problems}
Work related to the algebras $\Psi^\infty(\maG)$ used in the proof of our main 
result was done by Androulidakis, Connes, Debord, Mazzeo, Melrose, Monthubert, 
Ruzhansky, Schrohe, Schulze, Skandalis, and many others. See \cite{aln1, ConnesBook, DebordSkandalis, 
Mazzeo91, RuzhanskyBook, SchroheFrechet} for references. 
In that case, \m{\Psi^{-\infty}(\maG)} 
is the algebra generated by the conormal distributions of order \m{-\infty} 
on \m{A \to M}. In turn, \m{\Psi^{m}(\maG)}
is linearly generated by the conormal distributions of order \m{m} 
on \m{A \to M} and by \m{\Psi^{-\infty}(\maG).}
A generalization of these algebras in contained
in \cite{BCNQ}. In the case of manifolds with cylindrical ends, we have $\maV = \maV_b$
and $\Psi^{\infty}(\maG)$ consists of the properly supported operators 
in the $b$-calculus of Melrose and Schulze (and hence, it is a dense algebra 
in a suitable topology). Similarly, for asymptotically hyperbolic manifolds,
$\maV = \maV_0$ and $\Psi^{\infty}(\maG)$ consists of the properly supported 
operators in the $0$-calculus of Mazzeo and Schulze. Finally, for asymptotically 
euclidean manifolds, $\maV = \maV_{sc}$ and $\Psi^{\infty}(\maG)$ consists of 
the properly supported operators in the $SG$-calculus of Parenti and Schrohe,
the same calculus as the ``scattering calculus'' of Melrose.

\subsection{Problems}
We conclude by formulating two problems.

\begin{problem}[Connes, Bohlen-Schrohe] Find the index of a 
Fredholm operator \m{D \in \Psi^{m}(\maG).}
\end{problem}

\begin{problem}[Baldare-Benameur-Lesch-V.N., Ruzhansky] 
Let $G$ be a compact Lie group. 
Let \m{M} be compact smooth connected with a \m{G} action and \m{E} be a \m{G}-bundle. 
Let \m{\oM} be the Albin-Melrose compactification of the principal orbit
bundle (for the \m{G}) action and $\maG$ be its Lie groupoid. Describe \m{\Psi^m(M; E)^G} using 
\m{\Psi^{m}(\maG)} and the Ruzhansky calculus.
\end{problem}

\def\cprime{$'$}

\end{document}